\documentclass[10pt]{amsart}
\usepackage{amsmath,amssymb,latexsym,hyperref,url,graphicx,amsthm}

\newtheorem*{lemma}{Lemma}
\newtheorem*{proposition}{Proposition}
\newtheorem*{theorem}{Theorem}

\newcommand{\LP}{\textnormal{LP}}
\newcommand{\lp}{\textnormal{lp}}

\newcommand{\Eqn}[1]{Equation~(\ref{#1})}
\newcommand{\eqn}[1]{equation~(\ref{#1})}

\begin{document}

\title{Counting Interesting Elections}
\author{Lara K. Pudwell}
\author{Eric S. Rowland}
\date{October 30, 2009}

\begin{abstract}
We provide an elementary proof of a formula for the number of northeast lattice paths that lie in a certain region of the plane.  Equivalently, this formula counts the lattice points inside the Pitman--Stanley polytope of an $n$-tuple.
\end{abstract}

\maketitle \markboth{Lara Pudwell and Eric Rowland}{Counting Interesting Elections}

Suppose that on election day a TV news network of questionable morality wants to increase their viewership as polling results come in.  While the reporters cannot control the outcome of the election, they can control the order in which votes are reported to the public.  If one candidate is ahead in the tally throughout the entire day, viewership will wane since it is clear that she will win the election. On the other hand, a more riveting broadcast occurs when one candidate is ahead at certain times and the other candidate is ahead at others.  In fact, the network employs a group of psychologists and market analysts who have worked out certain margins they would like to achieve at certain points in the tally.  The director of programming needs to know the number of ways this can be done.

\section{The ballot problem}\label{ballotproblem}

We will work up to the general question by first examining the special (low ratings) case when one candidate has at least as many votes as the other throughout the tally.  This is the classical ``ballot problem'', in which candidate E and candidate N are competing for a public office.  Candidate E wins the election with $n$ votes.  How many ways are there to report the votes so that at all times during the tally N is not ahead of E?

We may represent the state of the tally at any moment by the pair $(x,y)$, where the coordinates $x$ and $y$ count the votes received by E and N respectively.  Then a tally consists of a sequence of points on the integer lattice in the plane made in steps of $E = \langle 1,0 \rangle$ and $N = \langle 0,1 \rangle$.  Such a sequence is called a \emph{northeast lattice path}.

\begin{figure}\label{restrictedpaths}
\begin{center}
\includegraphics{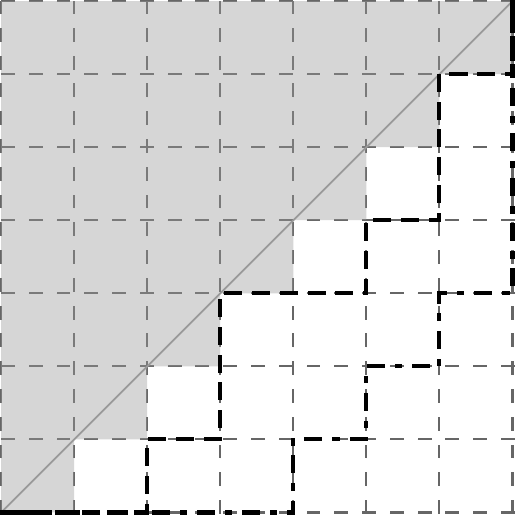}
\end{center}
\caption{Two northeast lattice paths from $(0,0)$ to $(7,7)$ restricted by $(EN)^7$.}
\end{figure}

We say that the lattice path $q$ is \emph{restricted} by the lattice
path $p$ if no part of $q$ lies directly above $p$.  For example,
Figure~\ref{restrictedpaths} shows two northeast lattice paths from
$(0,0)$ to $(n,n)$ that are restricted by the ``staircase'' $p =
ENEN \cdots EN$, or, equivalently, that do not go above the line $y
= x$. The ballot problem asks for the number $C_n$ of these paths.
(Note that if the tally ends at $(n,m)$, we may uniquely continue it
to a northeast lattice path ending at $(n,n)$.)

The ballot problem can be solved by constructing a simple recurrence.  Let $q$ be a northeast lattice path restricted by the staircase $p$.  Consider the point on $q$ where it first revisits the line $y=x$, and let $i$ be the $x$-coordinate of this point.  (This point exists since $q$ ends at $(n,n)$.)  For the upper path in Figure~\ref{restrictedpaths}, $i = 3$; for the lower path, $i = 7$.

Notice that since $q$ does not go above $y=x$ and begins at $(0,0)$
its first step is $E$; further, its last step before reaching the
point $(i,i)$ is $N$.  Therefore we may delete these steps to obtain
a northeast lattice path from $(1,0)$ to $(i,i-1)$ that does not go
above the line $y=x-1$.  There are $C_{i-1}$ ways to form such a
path, and there are $C_{n-i}$ ways to continue this path from
$(i,i)$ to $(n,n)$, so we have that $C_n = \Sigma_{i=1}^n C_{i-1}
C_{n-i}$. This we recognize as the familiar recurrence satisfied by
the Catalan numbers $C_n = \binom{2n}{n}/(n+1)$ \cite[Exercise
6.19(h)]{rS1999}, so we simply check that the initial condition $C_0
= 1$ agrees.

\section{Notation and theorem}\label{theorem}

We now consider a generalization of the ballot problem.  Let
$\LP(p)$ be the number of northeast lattice paths restricted by an
arbitrary northeast lattice path $p$ from $(0,0)$ to $(n,m)$.  The
path $p$ represents the network's predetermined restrictions on the
tally.  It was known by MacMahon \cite[p. 242]{macmahon} that the
sum of $\LP(p)$ over all such paths is
\[
    \sum_p \LP(p) = \frac{(m+n)!(m+n+1)!}{m!n!(m+1)!(n+1)!}.
\]
However, we are interested in computing $\LP(p)$ for specific $p$.

First we develop notation for lattice paths.

It is possible to represent a northeast lattice path as a word on $\{E, N\}$, such as
\[
    q = EENENNEENENNEN
\]
for the upper path in Figure~\ref{restrictedpaths}.  However, this representation is redundant, because the location of each $E$ step determines the path uniquely.

Therefore, we may represent a northeast lattice path by the sequence
of heights $q_i$ of the path along each interval from $x=i-1$ to
$x=i$.  For example, for the upper path in
Figure~\ref{restrictedpaths} we have $q = (0, 0, 1, 3, 3, 4, 6)$.
This representation is always a nondecreasing tuple of integers, and
it is our primary representation of lattice paths in this note.  A
lattice path $q = (q_1, q_2, \dots, q_n)$ is restricted by the
lattice path $p = (p_1, p_2, \dots, p_n)$ precisely when $q \leq p$
componentwise, i.e., $q_i \leq p_i$ whenever $1 \leq i \leq n$.

To write the main result, however, it turns out to be more natural to use still another representation of a northeast lattice path $p$ --- its difference sequence
\[
    \Delta p = (p_1, p_2 - p_1, \dots, p_n - p_{n-1}).
\]
Let $(v_1, v_2, \dots, v_n) = v = \Delta p$.  Since $p$ is a northeast lattice path, the entries of $v$ are nonnegative integers.  The entry $v_i$ is the number of $N$ steps taken along the line $x = i-1$, so we can think of this representation as determining a path by the location of each $N$ step.  The operator $\Delta$ has an inverse $\Sigma$, which produces the sequence of partial sums:
\[
    p = \Sigma v = (v_1, v_1 + v_2, \dots, v_1 + v_2 + \dots + v_n).
\]

The relationship between $p$ and $v = \Delta p$ can be interpreted in another way.  If $v = (v_1, v_2, \dots, v_n)$ is a tuple of nonnegative integers, the Pitman--Stanley polytope \cite{PS02} defined by $v$ is
\[
    \Pi_n(v) := \left\{ x \in \mathbb{R}_{\geq 0}^n : \text{$\Sigma x \leq \Sigma v$ componentwise} \right\}.
\]
Thus a tuple $x = (x_1, x_2, \dots, x_n)$ of nonnegative integers is a lattice point inside $\Pi_n(v)$ precisely when the northeast lattice path $\Sigma x$ is restricted by $\Sigma v$.  In other words, $\Delta$ provides a bijection from the northeast lattice paths restricted by $p$ to the lattice points in $\Pi_n(\Delta p)$.

We now return to the question at hand:  How many northeast lattice paths are restricted by the path $p = (p_1, p_2, \dots, p_{n-1}, p_n)$?  Equivalently, how many lattice points lie inside $\Pi_n(\Delta p)$?  One answer to this question is the following determinant enumeration.  Let $A = (a_{ij})$ be the $n \times n$ matrix with entries $a_{ij} = \binom{p_i + 1}{j - i + 1}$.  Then the number of northeast lattice paths restricted by $p$ is $\LP(p) = \det A$, as given by Kreweras \cite{kreweras} and Mohanty \cite[Theorem 2.1]{mohanty}.  This fact can be obtained from the triangular system of equations
\[
    \sum_{i=1}^{j+1} (-1)^{j-i+1} \binom{p_i + 1}{j - i + 1} \LP((p_1, p_2, \dots, p_{i-1})) =
    \begin{cases}
        1   & \text{if $j = 0$} \\
        0   & \text{if $j \geq 1$}
    \end{cases}
\]
for $0 \leq j \leq n$ (where $\LP(()) = 1$), which comes from an inclusion--exclusion argument; solve for $\LP((p_1, p_2, \dots, p_n))$ using Cramer's rule, and in the numerator expand by minors along the last column.

The following theorem presents a formula for $\LP(p)$ in which the lattice points in $\Pi_n((1, 1, \dots, 1, 1))$ play a central role.  This gives a non-determinantal formula for the number of northeast lattice paths restricted by $p$.  A generalization of the formula has been independently discovered by Gessel and by Pitman and Stanley \cite[Equation~(33)]{PS02} in more advanced contexts.  Our proof uses elementary combinatorial methods.

\begin{theorem}
Let $p$ be a northeast lattice path from $(0,0)$ to $(n,m)$, and let $v = \Delta p$.  The number of northeast lattice paths restricted by $p$ is
\begin{equation}\label{thm}
    \LP(p) = \sum_x \prod_{i=1}^n \binom{v_{n+1-i} + x_i - 1}{x_i},
\end{equation}
where the sum is over all $C_{n+1}$ lattice points $x$ in $\Pi_n((1, 1, \dots, 1, 1))$.
\end{theorem}

We immediately obtain two well-known results as special cases.  For $p = (m, m, \dots, m, m)$ we see that $v = \Delta p = (m, 0, \dots, 0, 0)$, which gives
\[
    \LP((m, m, \dots, m, m)) = \sum_x \binom{m+x_n-1}{x_n} \prod_{i=1}^{n-1} \binom{x_i-1}{x_i}.
\]
Since
\[
    \binom{x_i-1}{x_i} = \begin{cases}
        1 & \text{if $x_i = 0$} \\
        0 & \text{if $x_i \geq 1$}
    \end{cases}
\]
(from the generalization of the binomial theorem $(a + b)^m =
\sum_{j=0}^\infty \binom{m}{j} a^j b^{m-j}$ to $m=-1$), the only
nonzero terms in the sum come from lattice points of the form $(0,
0, \dots, 0, x_n)$, and therefore
\[
    \LP((m, m, \dots, m, m)) = \sum_{x_n=0}^n \binom{m+x_n-1}{x_n} = \binom{m+n}{n}
\]
as expected.

For $p = (1, 2, \dots, n-1, n)$ we recover the ballot problem.  Namely, $v = \Delta p = (1, 1, \dots, 1, 1)$, so
\[
    \LP((1, 2, \dots, n-1, n)) = \sum_x \prod_{i=1}^n \binom{x_i}{x_i} = \sum_x 1 = C_{n+1}.
\]

\Eqn{thm} allows one to compute $\LP(p)$ not only for explicit integer paths but for symbolic paths, and the resulting expressions have the pleasant property that they are written in the basis of rising factorials $a^{(m)} = a (a+1) \cdots (a+m-1)$.  For example, $\LP((v_1)) = v_1^{(0)} + v_1^{(1)} = 1 + v_1$.  For a general path of length $2$, we have
\begin{align*}
    \LP((v_1, v_1 + v_2)) &=
        v_2^{(0)} v_1^{(0)} +
        v_2^{(0)} v_1^{(1)} +
        \frac{1}{2} v_2^{(0)} v_1^{(2)} +
        v_2^{(1)} v_1^{(0)} +
        v_2^{(1)} v_1^{(1)} \\
    &= 1 + v_1 + \frac{1}{2} v_1 (v_1+1) + v_2 + v_2 v_1,
\end{align*}
and $\LP((v_1, v_1 + v_2, v_1 + v_2 + v_3))$ is
\begin{multline*}
    v_3^{(0)} v_2^{(0)} v_1^{(0)} +
    v_3^{(0)} v_2^{(0)} v_1^{(1)} \\
    + \frac{1}{2} v_3^{(0)} v_2^{(0)} v_1^{(2)} +
    \frac{1}{6} v_3^{(0)} v_2^{(0)} v_1^{(3)} +
    v_3^{(0)} v_2^{(1)} v_1^{(0)} +
    v_3^{(0)} v_2^{(1)} v_1^{(1)} +
    \frac{1}{2} v_3^{(0)} v_2^{(1)} v_1^{(2)} \\
    + \frac{1}{2} v_3^{(0)} v_2^{(2)} v_1^{(0)} +
    \frac{1}{2} v_3^{(0)} v_2^{(2)} v_1^{(1)} +
    v_3^{(1)} v_2^{(0)} v_1^{(0)} +
    v_3^{(1)} v_2^{(0)} v_1^{(1)} +
    \frac{1}{2} v_3^{(1)} v_2^{(0)} v_1^{(2)} \\
    + v_3^{(1)} v_2^{(1)} v_1^{(0)} +
    v_3^{(1)} v_2^{(1)} v_1^{(1)}.
\end{multline*}

Putting \eqn{thm} together with the determinantal formula for $\LP(p)$, we obtain a formula for a certain symbolic determinant in the same basis:
\[
    \det \binom{p_i + 1}{j-i+1}_{n \times n} = \sum_x \prod_{i=1}^n \frac{1}{x_i!} v_{n+1-i}^{(x_i)},
\]
where again $v = \Delta p$.

We note that Amdeberhan and Stanley \cite[Corollary 4.7]{amdeberhan-stanley} show that $\LP(p)$ also gives the number of monomials in the expanded form of the multivariate polynomial
\[
    \prod_{i=1}^n \sum_{j=1}^{p_i + 1} a_j
\]
in the variables $a_j$.  Moreover, $\LP(p)$ is the number of noncrossing matchings of a certain type \cite[Corollary 4.9]{amdeberhan-stanley}.

\section{Proof of the theorem}\label{proof}

Let $\lp(v)$ be the number of lattice points in $\Pi_n(v)$, where $v = (v_1, v_2, \dots, v_{n-1}, v_n)$.  That is, $\lp(v) = \LP(\Sigma v)$.  The following recurrence will be used.
\begin{proposition}
We have
\begin{equation*}
    \lp(v) = \begin{cases}
        1                               & \textnormal{if $n = 0$} \\
        \sum_{j=0}^{v_1} \lp((v_1+v_2-j, v_3, \dots, v_{n-1}, v_n)) & \textnormal{if $n \geq 1$.}
    \end{cases}
\end{equation*}
\end{proposition}

\begin{proof}
The only lattice point in $\Pi_0(())$ is $()$; hence $\lp(())=1$.

For $n \geq 1$, we partition the lattice points $w$ in $\Pi_n(v)$ according to the first entry $j = w_1$.  Since $w$ is a lattice point in $\Pi_n(v)$, then $w_1 + w_2 \leq v_1 + v_2$, so $w_2 \leq v_1 + v_2 - j$.  Therefore, lattice points $w = (j, w_2, \dots, w_{n-1}, w_n)$ in $\Pi_n(v)$ are in bijection (by deleting the first entry $j$) with lattice points in $\Pi_{n-1}((v_1+v_2-j, v_3, \dots, v_{n-1}, v_n))$.  Thus $\lp((v_1 + v_2 - j, v_3, \dots, v_{n-1}, v_n))$ is the number of lattice points in $\Pi_n(v)$ with first entry $j$, giving the recurrence.
\end{proof}

To prove the theorem, then, it suffices to show that \eqn{thm} satisfies this recurrence.  The base case $n = 0$ is easily checked, since the product is empty; we have
\[
    \sum_x \prod_{i=1}^n \binom{v_{n+1-i} + x_i - 1}{x_i} = \sum_x 1 = 1
\]
since again $\Pi_0(())$ has only one lattice point.

The remainder of this note is devoted to showing that for $n \geq 1$
\begin{multline}
    \sum_x \prod_{i=1}^n \binom{v_{n+1-i} + x_i - 1}{x_i} = \\
    \sum_{j=0}^{v_1} \sum_y \binom{v_1 + v_2 - j + y_{n-1} - 1}{y_{n-1}} \prod_{i=1}^{n-2} \binom{v_{n+1-i} + y_i - 1}{y_i},
\label{rec}
\end{multline}
where the left sum is over all $C_{n+1}$ lattice points $x$ in $\Pi_n((1, 1, \dots, 1, 1))$ and the right sum is over all $C_{n}$ lattice points $y$ in $\Pi_{n-1}((1, 1, \dots, 1))$.  We proceed by simplifying this equation until it becomes a statement about sums of binomial coefficients, given in the lemma below.

First interchange the two summations on the right side of \eqn{rec}.  Next, fix $y = (y_1, y_2, \dots, y_{n-1})$ on the right side, and break up the sum on the left according to the choice of $y$ in the following way.  The \emph{children} of $y = (y_1, y_2, \dots, y_{n-1})$ are the elements of the set
\[
    \{\, (y_1, y_2, \dots, y_{n-2}, y_{n-1}, 0) \,\} \cup \{\, (y_1, y_2, \dots, y_{n-2}, y_{n-1} - i, i + 1) : 0 \leq i \leq y_{n-1} \,\}.
\]
For example, the children of the lattice point $(0,3,2)$ are $(0,3,2,0)$, $(0,3,2,1)$, $(0,3,1,2)$, and $(0,3,0,3)$.  It is immediate that each lattice point $x$ has a unique parent $y$.

This definition is central to the proof.  The reason for defining
children in this way is that $x$ is a lattice point in $\Pi_{n}((1,
1, \dots, 1, 1))$ if and only if $x$'s parent is a lattice point in
$\Pi_{n-1}((1, 1, \dots, 1))$.  This property provides a many-to-one
correspondence between the $n$-dimensional lattice points in
$\Pi_n((1, 1, \dots, 1, 1))$ and the ($n-1$)-dimensional lattice
points in $\Pi_{n-1}((1, 1, \dots, 1))$. Using this correspondence
to break up \eqn{rec}, we obtain
\begin{multline}
    \sum_x \prod_{i=1}^n \binom{v_{n+1-i} + x_i - 1}{x_i} = \\
    \sum_{j=0}^{v_1} \binom{v_1+v_2-j+y_{n-1}-1}{y_{n-1}} \prod_{i=1}^{n-2} \binom{v_{n+1-i} + y_i - 1}{y_i}
\label{eq2}
\end{multline}
for each $y$, where the left sum is over all children $x = (x_1, x_2, \dots, x_{n-1}, x_n)$ of $y$.  It now suffices to prove \eqn{eq2} for a fixed $y$, since summing both sides of \eqn{eq2} over all $C_n$ lattice points $y$ in $\Pi_{n-1}((1, 1, \dots, 1))$ produces \eqn{rec}.

Note that if $x$ is a child of $y$ then $x_i = y_i$ for $1 \leq i \leq n-2$, so we may divide both sides of \eqn{eq2} by the product
\[
    \prod_{i=1}^{n-2} \binom{v_{n+1-i} + y_i - 1}{y_i}
\]
to obtain
\begin{equation}\label{eq3}
    \sum_x \binom{v_2 + x_{n-1} - 1}{x_{n-1}} \binom{v_1 + x_n - 1}{x_n} = \sum_{j=0}^{v_1} \binom{v_1 + v_2 - j + y_{n-1} - 1}{y_{n-1}}.
\end{equation}

We know what the children of $y$ look like, so the sum on the left side can be written as
\[
    \binom{v_2+y_{n-1}-1}{y_{n-1}} \binom{v_1+0-1}{0} + \sum_{i=0}^{y_{n-1}} \binom{v_2 + (y_{n-1}-i) - 1}{y_{n-1}-i} \binom{v_1 + (i+1) - 1}{i+1}.
\]
The first term in this expression, which corresponds to the child $(y_1, y_2, \dots, y_{n-1}, 0)$ of $y$, is equal to the $j = v_1$ term on the right side of \eqn{eq3}.  Removing this term from both sides leaves
\[
    \sum_{i=0}^{y_{n-1}} \binom{v_2+y_{n-1}-i-1}{y_{n-1}-i} \binom{v_1+i}{i+1} = \sum_{j=0}^{v_1-1} \binom{v_1+v_2-j+y_{n-1}-1}{y_{n-1}},
\]
which is proved in the following lemma under the substitution $a=v_1$, $b=v_2$, and $c=y_{n-1}$.

\begin{lemma}
Let $a$, $b$, and $c$ be nonnegative integers.  Then
\[
    \sum_{i=0}^c \binom{b+c-i-1}{c-i} \binom{a+i}{i+1} = \sum_{j=0}^{a-1} \binom{a+b+c-j-1}{c}.
\]
\end{lemma}

\begin{proof}
We show that both sides of the equation are equal to
\[
    \binom{a+b+c}{c+1} - \binom{b+c}{c+1}.
\]

The right side is a telescoping sum:
\begin{align*}
    \sum_{j=0}^{a-1} \binom{a+b+c-j-1}{c} &= \sum_{j=0}^{a-1} \left(\binom{a+b+c-j}{c+1} - \binom{a+b+c-j-1}{c+1}\right) \\
    &= \binom{a+b+c}{c+1} - \binom{b+c}{c+1}.
\end{align*}

The result for the left side follows from a generalization of the
Vandermonde identity, namely
\[
    \sum_{k=0}^f \binom{d+k}{k} \binom{e-k}{f-k} = \binom{d+e+1}{f}
\]
\cite[Problem 1.42(i)]{lovasz}.  The summand on the left side of
this equation counts the $(d+e+1-f)$-element subsets of
$\{1,2,\dots,d+e+1\}$ whose $(d+1)$st element is $d+k+1$ by choosing
$k$ of the first $d+k$ elements to be \emph{not} in the set and
$f-k$ of the last $e-k$ elements to be \emph{not} in the set.  The
right side counts all $(d+e+1-f)$-element subsets of
$\{1,2,\dots,d+e+1\}$ by selecting the elements \emph{not} in the
set.

Subtract $\binom{e}{f}$ from both sides of this equation and substitute $d=a-1$, $e=b+c$, $f=c+1$, and $k=i+1$ to obtain
\[
    \sum_{i=0}^c \binom{b+c-i-1}{c-i} \binom{a+i}{i+1} = \binom{a+b+c}{c+1} - \binom{b+c}{c+1}. \qedhere
\]
\end{proof}

Thus the director of programming may, for example, determine the likelihood that a random tally of votes will satisfy the network's needs.

\section*{Acknowledgements}

We thank Dimitrije Kostic for
generating our interest in this topic \cite{kostic}, Tewodros
Amdeberhan for pointing us to relevant literature, Doron Zeilberger
for encouraging us to simplify, and an anonymous referee for a
number of helpful suggestions.

\end{document}